\documentclass[11pt,oneside]{amsart}

\usepackage{amssymb}
\usepackage[latin1]{inputenc}
\usepackage[french]{babel}

\def\N{\mathbb{N}}
\def\Z{\mathbb{Z}}

\def\ZZ{\Z \oplus \Z}
\def\R{\mathbb{R}}
\def\H{\mathbb{H}}
\def\E{\mathbb{E}}
\def\KB{\mathbb{KB}}

\def\P{\pi_{1}}
\def\a{\alpha}

\def\l{\lambda}

\def\G{\Gamma}

\def\cal{\mathcal}

\def\ol{\overline}
\def\ul{\underline}

\def\wt{\widetilde}

\def\el{\'el\'ement}

\newtheorem{thm}{Théorème}[section]
\newtheorem{cor}{Corollaire}[section]
\newtheorem{lem}{Lemme}[section]

\title{Centralisateurs dans un groupe cristallographique 2-dimensionnel}

\begin{document}
\maketitle 
\begin{center}
{\sc Jean-Philippe PR\' EAUX}\footnote[1]{Centre de recherche de l'Armée de l'Air, Ecole de l'air, F-13661 Salon de
Provence air}\ \footnote[2]{Centre de Math\'ematiques et d'informatique, Universit\'e de Provence, 39 rue
F.Joliot-Curie, F-13453 marseille
cedex 13\\
\indent {\it E-mail :} \ preaux@cmi.univ-mrs.fr}
\end{center}
\begin{abstract}
Nous explicitons les centralisateurs dans un groupe discret cocompact d'isométrie du plan euclidien.
\end{abstract}

\section*{Introduction}
Nous donnons de façon explicite les centralisateurs dans un groupe cristallographique 2-dimensionnel (ou encore :
sous-groupes discrets cocompacts d'isométrie de $\E^2$, ou $\pi_1^{orb}$ d'orbiétés euclidiennes).

\section{Rappels sur les groupes Fuchsiens}

Soit $\G$ un groupe discret cocompact d'isométrie de $\E^2$, $\H^2$ ou $\mathbb{S}^2$. Il est connu que $\G$ admet une
des deux présentations suivantes où $p,q,g, \a_j\in \N^*$ et $\a_j>1$.\smallskip\\
Si $\G$ préserve l'orientation :
$$<{a}_{1},{b}_{1},\ldots,{a}_{g},{b}_{g},{c}_{1},
\ldots,{c}_{q},{d}_{1} ,\ldots,{d}_{p}\mid {c}_{j}^{\alpha _{j}}=1,\ (\prod_{i=1}^{g} \lbrack
{a}_{i},{b}_{i}\rbrack){c}_{1}\cdots {c}_{q} {d}_{1}\cdots {d}_{p}=1>$$
Si $\G$ renverse l'orientation :
$$<{a}_{1},\ldots ,{a}_{g},{c}_{1},\ldots,{c}_{q},{d}_{1}
,\ldots,{d}_{p}\mid {c}_{j}^{\alpha _{j}}=1,\ {a}_{1}^{2}\cdots {a}_{g}^{2}{c}_{1}\cdots {c}_{q}{d}_{1}\cdots {d}_{p}=1
>$$

 Ces groupes, parfois
appel\'es groupes Fuchsiens ou encore planaires, sont bien connus des analystes complexes, des g\'eom\`etres et des
topologues, et largement \'etudi\'es. Aussi, pour une introduction aux groupes Fuchsiens au sens où nous l'entendons
nous renvoyons le lecteur aux ouvrages de r\'ef\'erence  \cite{js}, \cite{hempel} et \cite{Lyndon}. Nous rappelerons
les r\'esultats qui nous serons utiles dans la
suite.\\

Topologiquement, le groupe Fuchsien $\G$ préservant (resp. renversant) l'orientation, est le groupe fondamental de
l'espace topologique $K$ que l'on obtient en retirant $q$ disques disjoints sur la surface $B$ orientable (resp. non
orientable) de genre $g$, et en recollant sur les $q$ composantes au bord de la surface obtenue, des disques, par des
applications de degr\'es respectifs $\a_1,\a_2,\ldots ,\a_q$.

Dans le cas o\`u $B$ est à  bord non vide, la derni\`ere relation peut \^ etre transform\'ee en
$${d}_{j}^{-1}={d}_{j+1}\cdots {d}_{p}.
({\prod_{i=1}^{g}} \lbrack {a}_{i},{b}_{i}\rbrack) .{c}_{1} \cdots {c}_{q}{d}_{1}\cdots {d}_{j-1}$$ ou
$${d}_{j}^{-1}={d}_{j+1}\cdots {d}_{p}.
({\prod_{i=1}^{g}}  {a}_{i}^2) .{c}_{1}\cdots {c}_{q}{d}_{1}\cdots {d}_{j-1}$$ pour tout $j=1\cdots n$. Un changement
de Tietze permet alors de supprimer cette relation et le g\'en\'erateur ${d}_{j}$. Ainsi $\G$ est le produit libre des
groupes cycliques engendr\'es par les g\'en\'erateurs ${a}_{i}$, ${b}_{j}$(dans le cas orientable), ${c}_{k}$, et
${d}_{l}$, \`a l'exception de l'un quelconque des ${d}_{l}$.

Il est bien connu qu'un produit libre de groupes cycliques contient un sous-groupe libre d'indice fini. En particulier,
si $\partial B\not= \varnothing$, alors $\G$ contient le groupe d'une surface \`a bord non vide, comme sous-groupe
d'indice fini. Ce r\'esultat reste vrai dans le cas o\`u la surface $B$ est fermée comme l'établit le th\'eor\`eme
suivant (th\'eor\`eme 12.2, \cite{hempel}).

\begin{thm}[Sous-groupe de surface d'indice fini]
\label{surf_fini} Soit $\G$ un groupe discret cocompact d'isométrie de $\E^2$, $\H^2$ ou $\mathbb{S}^2$. Si $\G$ admet
pour pr\'esentation :
$$<{a}_{1},{b}_{1},\ldots,{a}_{g},{b}_{g},{c}_{1},
\ldots,{c}_{q} \mid {c}_{j}^{\alpha _{j}}=1,\ (\prod_{i=1}^{g} \lbrack {a}_{i},{b}_{i}\rbrack){c}_{1} \cdots
{c}_{q}=1>\qquad\ (1)$$
ou
$$<{a}_{1},\ldots ,{a}_{g},{c}_{1},\ldots,{c}_{q}\mid
{c}_{j}^{\alpha _{j}}=1,\ {a}_{1}^{2}\cdots {a}_{g}^{2}{c}_{1}\cdots {c}_{q}=1 > \qquad\qquad\qquad (2)$$
alors, $\G$ contient le groupe d'une surface ferm\'ee $F$, comme sous-groupe d'indice fini, $\l$. De plus la
caract\'eristique d'Euler de $F$, est donn\'ee par la formule suivante :
$$\chi(F)=\begin{cases}
        \ \l (2-2g-\sum_{i=1}^{q}(1-\frac{1}{\a_i}))\qquad
&\text{dans le cas (1)}\\
        \ \l (2-g-\sum_{i=1}^{q}(1-\frac{1}{\a_i}))\qquad
&\text{dans le cas (2)}
        \end{cases}$$

\noindent En particulier, $\chi(F)>0$ ssi,\\
Cas (1) : $g=0$, et soit $q\leq 2$ soit $q=3$ et
$1/{\a_1}+1/{\a_2}+1/{\a_3}>1$.\\
Cas (2) : $g=1$, $q=1$ et $\a_1=2$ ou $g=0$  et soit $q\leq 2$ soit $q=3$ et $1/\a_1+1/\a_2+1/\a_3>1$.\smallskip\\
\noindent
et $\chi(F)=0$ seulement dans les cas suivants :\\
Cas (1) : $g=1$ et $q=0$ ; ou $g=0$, $q=3$, d'indices $(3,3,3),
(2,4,4)$ ou $(2,3,6)$ ; ou $g=0$, $q=4$, d'indices $(2,2,2,2)$.\\
Cas (2) : $g=2$ et $q=0$ ; ou $g=1$, $q=2$, d'indices $(2,2)$.
\end{thm}

\noindent {\bf Remarque 1 :} Si $\chi(F)\leq 0$, alors $\G$ est infini. Ainsi, si $B$ est ferm\'ee, $\G$ est fini ssi
$\chi(F)>0$. Si $\partial B\not= 0$, alors $\G$ est cyclique fini lorsque $g=0$, $q\geq 1$ et $p=1$, et infini sinon.\smallskip\\
 {\bf Remarque 2 :} Les
groupes finis obtenues (sous-groupes discrets d'isométrie de $S^2$, ou encore sous-groupes finis de $O(3)$)
sont :\\
-- $g=0$ et $q\leq 2$, $\G$ est cyclique fini.\\
-- $g=0$ et $q=3$, d'indices :\\
\indent -- $(2,2,n)$ ; $\G$ est le groupe diédral $D_{2n}$.\\
\indent -- $(2,3,3)$ ; $\G$ est le groupe tétraédral  $A_4$.\\
\indent -- $(2,3,4)$ ; $\G$ est le groupe octaédral $S_4$.\\
\indent -- $(2,3,5)$ ; $\G$ est le groupe icosaédral $A_5$.\\
-- $g=1$, $q=1$, $\a_1=2$ et $\G$ renverse l'orientation : $\G\cong\Z_4$.\\

 Dans un sens, un groupe Fuchsien est un "groupe de surface avec
torsion". Le théorème suivant \'etablit que si $\G$ est infini les seuls \el s de torsion sont \`a conjugaison pr\`es
les $c_i$. C'est la proposition 6.2 de \cite{Lyndon}, ou sous une forme plus r\'eduite la proposition II.3.6 de
\cite{js}.

\begin{thm}[Torsion dans un groupe Fuchsien infini]
\label{fuch_tor}
Soit $\G$, donn\'e par une des pr\'esentations (1) ou (2) figurant plus haut. Si $\G$ est infini, alors tout \el\ $u$
non trivial, de torsion,  est dans un conjugu\'e d'un des sous-groupes $<{c}_j>$, pour $j=1,\ldots ,q$. De plus il
s'\'ecrit de fa\c con unique sous la forme $u=v {c}_j^n v^{-1}$, avec $0<n<\a_j$, \emph{i.e.}, si il s'\'ecrit aussi
$u=w {c}_k^m w^{-1}$, avec $0<m<\a_k$, alors $j=k$, $n=m$, et $w=v{c}_j^t$, pour $t\in \Z$.
\end{thm}

\section{Centralisateurs des groupes Fuchsiens infinis}

Notons $\ol{B}$ la surface obtenue en retirant $q$ disques disjoints sur $B$. Si la surface $B$ est non orientable, il
en va de même de la surface $\ol{B}$. Consid\'erons le rev\^etement d'orientation de $\ol{B}$. Puisque tout lacet de
$\partial\ol{B}$ pr\'eserve l'orientation, ce rev\^etement s'\'etend \`a un rev\^etement \`a deux feuillets $\wt{K}$ de
$K$. Le groupe $\P(\wt{K})$ s'injecte dans $\P(K)=\G$, nous noterons $\ol{\G}$ son image. C'est un sous-groupe d'indice
2 dans $\G$.

Ce sous-groupe $\ol{\G}$ peut aussi se d\'efinir combinatoirement. Supposons encore que la surface $B$ soit non
orientable, et consid\'erons la pr\'esentation (2) de $\G$ donn\'ee auparavant. Un \el\ $g$ de $\G$ sera dit
${A}$-pair, si un mot le repr\'esentant contient un nombre pair d'occurences de lettres ${a}_1,{a}_1^{-1},\ldots
,{a}_g,{a}_g^{-1}$. Puisque les relateurs de $\G$ contiennent un nombre pair d'occurences de telles lettres, cette
notion est bien d\'efinie. Notons ${\ol{\G}}$ le sous-ensemble de $\G$ constitu\'e des \el s ${A}$-pairs ; il est
facile de v\'erifier que ${\Pi}$ est un sous-groupe d'indice 2 de $\G$. On retrouve la définition précédente : pour
voir celà, il suffit de remarquer que parmi les g\'en\'erateurs de $\G$ (donn\'es par la pr\'esentation (2)), seuls les
${a}_i$ ont pour repr\'esentants des lacets qui renversent l'orientation de $K$. \label{ulPi}  Si $B$ est orientable,
on pose ${\ol{\G}}=\G$.

Nous disposons maintenant du vocabulaire n\'ecessaire pour \'enoncer le r\'esultat caract\'erisant le centralisateur
d'un \el\ dans un groupe Fuchsien infini.

\begin{thm}[Centralisateurs dans un groupe Fuchsien infini]
 \label{fuch_cent}
 Soit $\G$ un groupe Fuchsien infini, et $u$ un \el\ non trivial de $\G$.
Si $u$ est d'ordre infini, alors son centralisateur $\cal{Z}(u)$ est soit isomorphe au groupe de la bouteille de klein,
soit libre ab\'elien de rang $1$ ou $2$. Si de plus $u\not\in {\ol{\G}}$, alors $\cal{Z}(u)$ est cyclique infini.

Si $u$ est d'ordre fini, alors son centralisateur est cyclique fini. Plus pr\'ecis\'ement (cf. proposition
\ref{fuch_tor}), pour $i=1,\ldots ,q$, et $0<n<\a_i$, $\cal{Z}({c}_i^n)=<{c}_i>$.
\end{thm}

\noindent \textbf{Remarque :} Si $\G$ est fini, le centralisateur d'un \el\ n'est pas n\'ecessai\-rement cyclique.
Ainsi par exemple si $B=S^2$, $q=3$ d'indices $(2,2,2)$, $\G\cong \Z_2\times Z_2$ est abélien non cyclique, et $\G$ est
le centralisateur de chacun de ses \el s.\\

\noindent \textsl{D\'emonstration} Notons $\ol{\G}$ le sous-groupe canonique de $\G$. Consid\'erons le centralisateur
$\cal{Z}(u)$ d'un \el\ $u$ de $\G$.

Si $u$ est d'ordre infini  $\cal{Z}(u)$ est un groupe infini ayant un centre non trivial.
Maintenant, d'une part tout sous-groupe d'un groupe Fuchsien est un groupe Fuchsien (cf. \S 3 de \cite{js}), et ainsi
$\cal{Z}(u)$ est un groupe Fuchsien. D'autre part, un groupe Fuchsien infini ayant un centre non trivial est soit le
groupe de la bouteille de Klein, soit libre ab\'elien de rang 1 ou 2. Si son centre  n'est pas contenu dans le
sous-groupe canonique $\cal Z(u)\cap\ol{\G}$ de $\cal Z(u)$  alors $G$ est cyclique infini (cf. proposition II.3.11,
\cite{js}). Ainsi si $u\not\in \ol{\G}$ alors $\cal{Z}(u)$ est cyclique infini. Ceci d\'emontre la premi\`ere partie de
la proposition.

Si $u$ est d'ordre fini, alors $\cal{Z}(u)$ ne peut pas contenir d'\el\ d'ordre infini. En effet, dans le cas contraire
$\cal{Z}(u)$ contiendrait le sous-groupe $\Z\times \Z_n$, qui n'est pas un groupe Fuchsien (lemme II.3.10 \cite{js}),
ce qui est contradictoire. Avec le théorème \ref{fuch_tor}, si $u$ est d'ordre fini, $u$ est conjugu\'e  \`a ${c}_i^n$,
pour $1\leq i\leq q$ et $n\in \Z_\ast$, et nous  supposerons sans perte de généralité que $u={c}_i^n$. Consid\'erons un
\el\ non trivial $v\in \cal{Z}(u)$, il est de torsion, et donc avec le théorème \ref{fuch_tor}, $v=a{c}_j^ma^{-1}$.
Puisque $u$ et $v$ commutent, $u={c}_i^n= a{c}_j^ma^{-1} {c}_i^n a{c}_j^{-m}a^{-1}$ et donc avec le théorème
\ref{fuch_tor} il existe $k\not= 0$, $v=a{c}_j^ma^{-1}={c}_i^k$
et donc $v \in <{c}_i>$. Ainsi, pour $n\not=0$, $\cal{Z}({c}_i^n)=<{c}_i>$.\hfill$\square$\\

Soit $\G$ un groupe discret cocompact d'isométrie de $\E^2$ ou $\H^2$. Soit $\G$ contient $\ZZ$ soit $\G$ est un
sous-groupe discret de $PSL(2,\R)=Isom^+(\H^2)$. Dans le dernier cas, $\G$ est soit produit libre de groupes cycliques
soit sans torsion, et dans tous les cas hyperbolique au sens de Gromov. Le centralisateur d'un \el\ est cyclique, et
l'on peut algorithmiquement le déterminer pour un \el\ arbitraire en appliquant l'algorithme figurant en \cite{moi}.

Les cas restant correspondent au cas euclidien, ou encore  \`a $B$ ferm\'ee, et $\G$ contenant un sous-groupe d'indice
fini, isomorphe au groupe d'une surface $F$ avec $\chi(F)=0$. Il existe 7 groupes Fuchsiens v\'erifiant ces conditions,
ils sont caract\'eris\'es dans le th\'eor\`eme \ref{surf_fini}. Dans ce cas, $\G$ contient $\ZZ$ comme sous-groupe
d'indice fini et le centralisateur d'un \el\ n'est plus n\'ecessairement cyclique. Le lemme suivant fournit une
pr\'esentation plus ad\'equate pour chacun de ces groupes, en ceci qu'elle explicite le sous-groupe libre ab\'elien de
rang 2 d'indice fini, et une d\'ecomposition en un produit semi-direct d'un groupe de surface euclidienne par un groupe
fini.

\begin{lem}[Pr\'esentation pour un groupe Fuchsien
cristallographique] \label{crist_pres}  Consid\'erons les groupes Fuchsiens suivants :
\begin{gather*}
G_1=<{a}_1,{a}_2\mid {a}_1^2{a}_2^2=1>\\
G_2=<{c}_1,{c}_2,{c}_3\mid
{c}_1^2={c}_2^2={c}_3^2=({c}_1{c}_2{c}_3)^2=1>\\
G_3=<{c}_1,{c}_2\mid {c}_1^3={c}_2^3=({c}_1{c}_2)^3=1>\\
G_4=<{c}_1,{c}_2,\mid {c}_1^4={c}_2^4=({c}_1{c}_2)^2=1>\\
G_5=<{c}_1,{c}_2,\mid {c}_1^6={c}_2^3=({c}_1{c}_2)^2=1>\\
G_6=<{a},{c}, \mid {c}^2=({a}^2{c})^2=1>
\end{gather*}
Ils admettent aussi pour pr\'esentations :
\begin{gather*}
G_1\cong <{a},t_1,t_2\mid {a}^2=t_1,\ [t_1,t_2]=1,\
{a}t_2{a}^{-1}=t_2^{-1}>\\
G_2\cong <c,t_1,t_2\mid c^2=1,\ [t_1,t_2]=1,\ ct_1c^{-1}=t_1^{-1},\
ct_2c^{-1}=t_2^{-1}>\\
G_3\cong <c,t_1,t_2\mid c^3=1,\ [t_1,t_2]=1,\ ct_1c^{-1}=t_2,\
ct_2c^{-1}=t_1^{-1}t_2^{-1}>\\
G_4\cong <c,t_1,t_2\mid c^4=1,\ [t_1,t_2]=1,\ ct_1c^{-1}=t_2,\
ct_2c^{-1}=t_1^{-1}>\\
G_5\cong <c,t_1,t_2\mid c^6=1,\ [t_1,t_2]=1,\ ct_1c^{-1}=t_2,\
ct_2c^{-1}=t_1^{-1}t_2>\\
G_6\cong <{a}, {c},t_1,t_2\mid\
 {a}^2=t_1,\ [t_1,t_2]=1,\ {a}t_2{a}^{-1}=t_2^{-1},
\qquad\qquad\qquad\qquad\\
\qquad\qquad\qquad\qquad c^2=1,  ct_1c^{-1}=t_1^{-1},\ ct_2c^{-1}=t_2^{-1}, c\,{a}\,c^{-1}=t_2{a}^{-1}>
\end{gather*}
\end{lem}

\noindent \textbf{Remarques :} -- Les pr\'esentations obtenues explicitent le sous-groupe $\ZZ$ d'indice fini et la
structure de produit semi-direct pour chacun de ces groupes.

-- Bien s\^ur $G_1$ est le groupe de la bouteille de Klein,  $G_2,G_3,G_4,G_5$ sont des produits semi-directs de $\ZZ$,
respectivement par $\Z_2$, $\Z_3$, $\Z_4$, $\Z_6$. Quant \`a $G_6$, c'est un produit semi-direct du groupe de la
bouteille de Klein par $\Z_2$ mais aussi une extension finie de $\ZZ$ par $\Z_2\oplus \Z_2$.

-- Les sous-groupes d'orientation $\ol{G}_1,\ol{G}_6$ de $G_1$ et $G_6$ sont engendr\'es respectivement
par ${a}^2,t_1,t_2$ et ${a}^2,c,t_1,t_2$.\\

\noindent\textsl{D\'emonstration.} Proc\'eder par changements de Tietze. Nous n'indiquons que comment exprimer les
nouveaux générateurs en fonction des anciens.

\noindent{Pour} ${G_1.}$
 Poser $a=a_1$, $t_1={a}_1^2$, et $t_2={a}_1^{-1}{a}_2^{-1}$.

\noindent{Pour} ${G_2.}$ Posons $c={c}_1$, $t_1={c}_2{c}_1$, et $t_2={c}_2{c}_3$.

\noindent{Pour} ${G_3.}$
 Posons $c={c}_1^{-1}$, $t_1={c}_1^{-1}{c}_2$,
et $t_2={c}_1{c}_2{c}_1$.

\noindent {Pour} ${G_4.}$ Posons $c={c}_1$, $t_1={c}_2{c}_1^{-1}$, et $t_2={c}_1{c}_2{c}_1^2$.

\noindent
 {Pour}
${G_5.}$ Posons $c={c}_1$, $t_1={c}_1^{-2}{c}_2$, et $t_2={c}_1^{-1}{c}_2{c}_1^{-1}$.

 \noindent{Pour}
${G_6.}$ Posons  $t_1={a}^2$, et $t_2=({c}{a})^2$. \hfill$\square$

Nous pouvons dès-lors expliciter les centralisateurs dans ces groupes. Une preuve calculatoire peut se trouver dans ma
thèse de doctorat (\cite{thesis}).

\begin{thm}[Centralisateurs des groupes Fuchsiens cristallographiques]
\label{cent_crist} Consid\'erons  les groupes $G_1,G_2,G_3,G_4,G_5$ et $G_6$ apparaissant dans le lemme
\ref{crist_pres}, et munissons-les des  pr\'esentations qu'il nous fournit. Elles explicitent pour chacun des $G_i$ un
sous-groupe d'indice fini, libre ab\'elien de rang 2, engendr\'e par $t_1,t_2$; notons le $H_i$. Nous d\'ecrivons le
centralisateur d'un \el\ quelconque de ces groupes.
 Si $u\in G_i$, notons
$\cal{Z}(u)$ son centralisateur dans $G_i$.\\

\noindent $\ul{Pour\ G_1}$: Tout \el\ de $G_1$ s'\'ecrit de fa\c con unique, sous la forme $t_1^{n_1}t_2^{n_2}$ ou
$t_1^{n_1}t_2^{n_2}{a}$. Le centralisateur d'un \el\ non trivial $u\in G_1$,
est d\'ecrit par :\\
\begin{center}
\begin{tabular}{|c||c|}
\hline $\left.\begin{array}{ll}
    u\in H_1\\
    u=t_1^{n_1}t_2^{n_2}
    \end{array}\right.$
& $\cal{Z}(u)=\left\lbrace
\begin{array}{ll}
    G_1\cong\P(\KB_2)\qquad\ si\ u\in<t_1>\\
    H_1\cong \ZZ\qquad \quad\ si\ u\not\in <t_1>
\end{array}\right.$
\\
\hline $\left.\begin{array}{ll}
    u\in H_1.{a}\\
    u=t_1^{n_1}t_2^{n_2}{a}
    \end{array}\right.$
& $\cal{Z}(u)=<t_2^{n_2}{a}>$\\
\hline
\end{tabular}
\end{center}

\noindent$\ul{Pour\ G_2}$: Tout \el\ de $G_2$ s'\'ecrit de fa\c con unique sous la forme $t_1^{n_1}t_2^{n_2}$ ou
$t_1^{n_1}t_2^{n_2}c$. Le centralisateur d'un \el\ non trivial $u$
de $G_2$ est d\'ecrit par :\\
\begin{center}
\begin{tabular}{|c||c|}
\hline $\left.\begin{array}{ll}
    u\in H_2\\
    u=t_1^{n_1}t_2^{n_2}
    \end{array}\right.$
& $\cal{Z}(u)=H_2\cong\ZZ$
\\
\hline $\left.\begin{array}{ll}
    u\in H_2.c\\
    u=t_1^{n_1}t_2^{n_2}c
    \end{array}\right.$
&
$\cal{Z}(u)=<u>$\\
\hline
\end{tabular}
\end{center}

\noindent$\ul{Pour\ G_3}$: Tout \el\ de $G_3$ s'\'ecrit de fa\c con unique sous la forme $t_1^{n_1}t_2^{n_2}$, ou
$t_1^{n_1}t_2^{n_2}c$, ou $t_1^{n_1}t_2^{n_2}c^2$. Le centralisateur d'un \el\ non trivial $u$ de $G_3$ est d\'ecrit
par :
\begin{center}
\begin{tabular}{|c||c|}
\hline $\left.\begin{array}{ll}
    u\in H_3\\
    u=t_1^{n_1}t_2^{n_2}
    \end{array}\right.$
&
$\cal{Z}(u)=H_3\cong\ZZ$\\
\hline $\left.\begin{array}{ll}
    u\in H_3.c\\
    u=t_1^{n_1}t_2^{n_2}c
    \end{array}\right.$
&
$\cal{Z}(u)=<u>$\\
\hline $\left.\begin{array}{ll}
    u\in H_3.c^2\\
    u=t_1^{n_1}t_2^{n_2}c^2
    \end{array}\right.$
&
$\cal{Z}(u)=<t_1^{n_2}t_2^{n_2-n_1}c>$\\
\hline
\end{tabular}
\end{center}

\noindent$\ul{Pour\ G_4}$: Tout \el\ de $G_4$ s'\'ecrit de fa\c con unique sous la forme $t_1^{n_1}t_2^{n_2}$, ou
$t_1^{n_1}t_2^{n_2}c$, ou $t_1^{n_1}t_2^{n_2}c^2$, ou $t_1^{n_1}t_2^{n_2}c^3$. Le centralisateur d'un \el\ non trivial
$u$ de $G_4$ est d\'ecrit par :
\begin{center}
\begin{tabular}{|c||c|}
\hline $\left.\begin{array}{ll}
    u\in H_4\\
    u=t_1^{n_1}t_2^{n_2}
    \end{array}\right.$
& $\cal{Z}(u)=H_4\cong \ZZ$
\\
\hline $\left.\begin{array}{ll}
    u\in H_4.c\\
    u=t_1^{n_1}t_2^{n_2}c
    \end{array}\right.$
&
$\cal{Z}(u)=<u>$\\
\hline $\left.\begin{array}{ll}
    u\in H_4.c^2\\
    u=t_1^{n_1}t_2^{n_2}c^2
    \end{array}\right.$
& $\cal{Z}(u)=\left\lbrace
    \begin{array}{ll}
    <t_1^{\frac{n_1+n_2}{2}}
t_2^{\frac{n_2-n_1}{2}}c>  \qquad\quad\, \text{si $n_1+n_2$ pair}\\
    <u> \qquad\qquad\qquad\qquad \text{sinon}
    \end{array}\right.$
\\
\hline $\left.\begin{array}{ll}
    u\in H_4.c^3\\
    u=t_1^{n_1}t_2^{n_2}c^3
    \end{array}\right.$
&
$\cal{Z}(u)=<t_1^{n_2}t_2^{-n_1}c>$\\
\hline
\end{tabular}
\end{center}

\noindent$\ul{Pour\ G_5}$: Tout \el\ de $G_5$ s'\'ecrit de fa\c con unique sous la forme $t_1^{n_1}t_2^{n_2}$, ou
$t_1^{n_1}t_2^{n_2}c$, ou $t_1^{n_1}t_2^{n_2}c^2$, ou $t_1^{n_1}t_2^{n_2}c^3$, ou $t_1^{n_1}t_2^{n_2}c^4$, ou
$t_1^{n_1}t_2^{n_2}c^5$. Le centralisateur d'un \el\ non trivial $u$ de $G_5$ est d\'ecrit par :
\begin{center}
\begin{tabular}{|c||c|}
\hline $\left.\begin{array}{ll}
    u\in H_5\\
    u=t_1^{n_1}t_2^{n_2}
    \end{array}\right.$
& $\cal{Z}(u)=H_5\cong \ZZ$
\\
\hline $\left.\begin{array}{ll}
    u\in H_5.c\\
    u=t_1^{n_1}t_2^{n_2}c
    \end{array}\right.$
&
$\cal{Z}(u)=<u>$\\
\hline $\left.\begin{array}{ll}
    u\in H_5.c^2\\
    u=t_1^{n_1}t_2^{n_2}c^2
    \end{array}\right.$
& $\cal{Z}(u)=\left\lbrace
    \begin{array}{ll}
    <t_1^{\frac{2n_1+n_2}{3}}
t_2^{\frac{n_2-n_1}{3}}c>\quad\ si\ 3/(n_1-n_2)\\
    <u>\qquad \qquad \qquad \qquad \ sinon
    \end{array}\right.$
\\
\hline $\left.\begin{array}{ll}
    u\in H_5.c^3\\
    u=t_1^{n_1}t_2^{n_2}c^3
    \end{array}\right.$
& $\cal{Z}(u)=\left\lbrace
    \begin{array}{ll}
    <t_1^{\frac{n_1+n_2}{2}}
t_2^{-\frac{n_1}{2}}c>\quad si\ n_1\ et\ n_2\ pairs\\
    <u>\qquad \qquad \qquad \qquad \ sinon
    \end{array}\right.$
\\
\hline $\left.\begin{array}{ll}
    u\in H_5.c^4\\
    u=t_1^{n_1}t_2^{n_2}c^4
    \end{array}\right.$
& $\cal{Z}(u)=\left\lbrace
    \begin{array}{ll}
    <t_1^{\frac{n_1+2n_2}{3}}
t_2^{-\frac{2n_1+n_2}{3}}c>\quad si\ 3/(n_1-n_2)\\
    <t_1^{n_1+n_2}t_2^{-n_1}c^2>\qquad \qquad \, sinon
    \end{array}\right.$
\\
\hline $\left.\begin{array}{ll}
    u\in H_5.c^5\\
    u=t_1^{n_1}t_2^{n_2}c^5
    \end{array}\right.$
& $\cal{Z}(u)=<t_1^{n_2}t_2^{-(n_1+n_2)}c>$
\\
\hline
\end{tabular}
\end{center}

\noindent$\ul{Pour\ G_6}$: Tout \el\ de $G_6$ s'\'ecrit uniquement sous la forme $t_1^{n_1}t_2^{n_2}$ ou
$t_1^{n_1}t_2^{n_2}{a}$, ou $t_1^{n_1}t_2^{n_2}{c}$, ou $t_1^{n_1}t_2^{n_2}{a}\,{c}$. Notons $K_6$ le sous-groupe de
$G_6$ engendr\'e par ${a},t_1,t_2$; il est isomorphe au groupe de la bouteille de Klein. Le centralisateur d'un \el\
non trivial $u$ de $G_6$ est d\'ecrit par :
\begin{center}
\begin{tabular}{|c||c|}
\hline $\left.\begin{array}{ll}
u\in H_6\\
u=t_1^{n_1}t_2^{n_2}
\end{array}\right.$
& $\cal{Z}(u)=\left\lbrace
    \begin{array}{lll}
    <{a}c\, ,t_1>\cong \P(\KB_2)
\qquad si\ u\in <t_2>\\
    K_6\cong \P(\KB_2) \qquad \qquad \ \ si\ u\in <t_1>\\
    H_6\cong \ZZ \qquad \qquad \qquad \ sinon
    \end{array}\right.$
\\
\hline $\left.\begin{array}{ll}
u\in H_6.{a}\\
u=t_1^{n_1}t_2^{n_2}{a}
\end{array}\right.$
& $\cal{Z}(u)=<t_2^{n_2}{a}>$
\\
\hline $\left.\begin{array}{ll}
u\in H_6.c\\
u=t_1^{n_1}t_2^{n_2}{c}
\end{array}\right.$
& $\cal{Z}(u)=<u>$
\\
\hline $\left.\begin{array}{ll}
u\in H_6.{a}{c}\\
u=t_1^{n_1}t_2^{n_2}{a}{c}
\end{array}\right.$
& $\cal{Z}(u)=<t_1^{n_1}{a}{c}>$
\\
\hline
\end{tabular}
\end{center}
\end{thm}

\vskip 0.2cm

\begin{cor}
Sous les mêmes hypothèses que précédemment, $G_1$ a un centre cyclique infini engendré par $t_1$, tandis que $G_2$,
$G_3$, $G_4$, $G_5$, $G_6$ ont un centre trivial.
\end{cor}


\vskip 0.4cm


\end{document}